\documentclass[11pt]{article}

\usepackage{amsmath,amsfonts,amssymb,amsthm,bbm,latexsym,mathrsfs}
\usepackage{graphicx,color,epsfig,fancyhdr,dsfont, ulem}
\usepackage{enumerate}
\usepackage{hyperref}
\usepackage{indentfirst}
\usepackage{amsmath,amscd}
\usepackage{caption}
\usepackage{graphicx, subfig}
\usepackage[all]{xy}
\usepackage[dvipsnames]{xcolor}
\title{On minimal flows and definable amenability in some distal $\NIP$ theories}
\date{\today}

\author{Ningyuan Yao \\
  \text{Fudan University} \and Zhentao Zhang\\
  \text{Fudan University}}
\date{}

\newtheorem{theorem}{Theorem}[section]
\newtheorem{prop}[theorem]{Proposition}
\newtheorem{definition}[theorem]{Definition}
\newtheorem{rmk}[theorem]{Remark}
\newtheorem{lemma}[theorem]{Lemma}
\newtheorem{coro}[theorem]{Corollary}
\newtheorem{fact}[theorem]{Fact}

\newtheorem{Conj}{Conjecture}

\newcommand{\tp}{\mathrm{tp}}
\newcommand{\Th}{\mathrm{Th}}
\newcommand{\cl}{\mathrm{cl}}
\newcommand{\Gen}{\mathrm{Gen}}
\newcommand{\ext}{\mathrm{ext}}
\newcommand{\dcl}{\mathrm{dcl}}
\newcommand{\Aut}{\mathrm{Aut}}
\newcommand{\Q}{\mathbb{Q}_p}
\newcommand{\Z}{\mathbb{Z}}
\newcommand{\M}{\mathbb{M}}

\newcommand{\N}{\mathbb{N}}

\newcommand{\Ga}{\mathbb{G}_\text{a}}
\newcommand{\Gm}{\mathbb{G}_\text{m}}

\newcommand{\pCF}{p\mathrm{CF}}

\newcommand{\Stab}{\mathrm{Stab}}
\newcommand{\WG}{\mathrm{WG}}
\newcommand{\AP}{\mathrm{AP}}
\newcommand{\SL}{\mathrm{SL}}
\newcommand{\SO}{\mathrm{SO}}
\newcommand{\st}{\mathrm{st}}
\newcommand{\I}{\mathcal{I}}
\newcommand{\J}{\mathcal{J}}
\newcommand{\RCF}{\mathrm{RCF}}
\newcommand{\dfg}{\mathrm{dfg}}
\newcommand{\fsg}{\mathrm{fsg}}
\newcommand{\NIP}{\mathrm{NIP}}

\newcommand{\sq}{\subseteq}

\begin{document}

\maketitle

\begin{abstract}

We study the  definable topological dynamics $(G(M), S_G(M))$ of a definable group acting on its type space, where $M$ is either an $o$-minimal structure or a $p$-adically closed field, and $G$ a  definably amenable group. We  focus on the problem raised in \cite{Newelski-I} of whether weakly generic types coincide with almost periodic types, showing that the answer is positive when $G$ has boundedly many global weakly generic types. We also give two ``minimal counterexamples'' where $G$ has unboundedly many global weakly generic types,  extending  the main results in \cite{Pillay-Yao-mini-flow} to a more general context.
\end{abstract}

\section{Introduction}

In model theory, we study a group $G$ definable in a structure $M$ and the action of $G$ on its type space $S_G(M)$, which is the collection of all types over $M$ containing the formula defining $G$. The space of generic types, introduced by Poizat as a generalization of the notion of generic points in an algebraic group, plays a heart role when $\Th(M)$ is stable. But  for the unstable cases, the generic types may not exist. So
various of weakenings of the generic were introduced to the unstable environment to generalize the properties of stable groups to  the unstable context. The notion of weakly generic types introduced by Newelski in \cite{Newelski-I}, which exists in any context, is a suitable substitution for generic types. We say that a definable set $X$ is \emph{weakly generic} if there is a non-generic definable set $Y$ such that $X\cup Y$ is generic, where a definable set is \emph{generic} if finitely many its translates cover the whole group. Newelski studied the action of a definable group on its type space in the  topological dynamics point of view, and tried to link the invariants suggested by topological dynamics with model-theoretic invariants.

Almost periodic types are one of the  “new” objects suggested by topological dynamics. We say that a type $p\in S_G(M)$ is \emph{almost periodic} if the closure of its $G(M)$-orbit is a minimal subflow of $S_G(M)$. Newelski proved that the space of almost periodic types coincides with the closure of space of weakly generic types, and when the generic types exist, almost periodics coincides with the weakly generics (see  Corollary 1.9 and Remark 1.10 in \cite{Newelski-I}). An example was given where the two classes differ and the problem was
explicitly raised (Problem 5.4 of \cite{Newelski-I}) of finding an $o$-minimal or even just $\NIP$ example. Newelski's question is restated in \cite{CS-amenable-$NIP$-group} (Question 3.35) in the special case of definably amenable groups in $\NIP$ theories.

When $M$ is an $o$-minimal expansion of a real closed field and $G$ is a definably amenable group definable over $M$, Pillay and Yao proved in \cite{Pillay-Yao-mini-flow} that weakly generics coincide with  almost periodics when the torsion free part of $G$ has dimension one.  They also construct a counter-example when $G=S^1\times (\mathbb R,+)^2$, to show that the set of weakly generic types properly contains the set of almost periodic types. The existence of a non-stationary weakly generic type in $({\mathbb R}^2,+)$ plays a crucial role in the construction of the counter-example, where a weakly generic type is by definition stationary if it has a unique global weakly generic extension (see \cite{Petri}, Definition 3.41).

The current paper continues this line of work,  extending the main results of \cite{Pillay-Yao-mini-flow} to a rather broader context, where $T$ is  $\pCF$ of $p$-adically closed fields or a complete $o$-minimal expansion of the theory $\RCF$ of real closed fields.  Let $M$ be a model of $T$ and $M=\Q$ when $T=\pCF$. We consider the case that $G$ is a $M$-definable group admitting a $M$-definable short exact sequence
\begin{equation}\label{equ-dfg-fsg-decom}
    1\rightarrow H\rightarrow G\rightarrow C\rightarrow 1,
\end{equation}
where $C$ is definably compact (see \cite{PS-definably-comp} and \cite{OP-comp-p-adic} for definitions) and $H$ is ``totally non-compact''. We say that $H$  is \emph{totally non-compact} if there is a definable normal sequence
\begin{equation}\label{tot-non-cpt}
    H_0\lhd H_1\lhd ...\lhd H_i\lhd...\lhd H_n=H
\end{equation}
such that $H_0$ is finite and  each $H_{i+1}/H_{i}$ is a one-dimensional non-definably compact group. Note that any group definable over $M$ admitting a decomposition as in (\ref{equ-dfg-fsg-decom}) is definably amenable (see Remark \ref{DA-by-DA}).

If $H$ is defined in an $o$-minimal expansion of $\RCF$, then $H$ is torsion-free iff it is totally non-compact and each $H_i$ as in \ref{tot-non-cpt} is definably connected (see \cite{Edmu}).
Any group $G$ definable in an $o$-minimal expansion of $\RCF$ is definably amenable  iff it has a decomposition as in (\ref{equ-dfg-fsg-decom}) (see \cite{CP-connected-component}).

When $H$ is definable over $\Q$, $H$  is eventually trigonalizable  over $\Q$ \cite{PY-dfg-groups}. A recent result of \cite{JY-abelian-pcf-group} shows that any $\Q$-definable abelian group  also has  a decomposition as in (\ref{equ-dfg-fsg-decom}).

We can also describe the ``totally non-compact'' and  ``definably compact'' via the model-theoretic invariants of ``$\dfg$'' and ``$\fsg$'' respectively, in both $o$-minimal and $\pCF$ contexts.  Recall that a group $G$ has finitely satisfiable generics ($\fsg$) or definable $f$-generics ($\dfg$) if there is a global type $p$ on $G$ and a small model $M$ such that every left translate of $p$ is finitely satisfiable in $M$ or definable over $M$, respectively. Suppose that $M$ is either $\Q$ or an $o$-minimal expansion of a real closed field, and $H$ is definable over $M$, then $H$ is totally non-compact iff $H$ has $\dfg$ \cite{CP-connected-component,PY-dfg-groups}, and $C$ is definably compact iff $C$ has $\fsg$ \cite{HPP-$NIP$-JAMS, OP-comp-p-adic}. So we assume in this paper that $G$ admits a ``$\dfg$-$\fsg$'' short exact sequence
\begin{equation}\label{equ-fsg-fsg-1}
     1\rightarrow H\rightarrow G\rightarrow C\rightarrow 1,
\end{equation}
namely, $G$ is an extension of a $\fsg$ group $C$ by a $\dfg$ group $H$.

The structure of this paper is analogous to \cite{Pillay-Yao-mini-flow}. We study the topological dynamics in both global and local environments. What we call the global context is where $\M$ is a monster model of $T$. What we call the local context is where $M$ is any model of $T$, and we pass to the Shelah expansion $M_0=M^\ext$ of $M$ by the externally definable sets and consider the action of $G(M_0)$ on $S_G(M_0)$.

We introduce a new terminology first. As it was mentioned above a weakly generic type  over $A$ is stationary if it has just one extension to a complete weakly generic type over the monster model $\M\supset A$. In this paper, we call a definable group $X$ is \emph{stationary} if every weakly generic type over every small model $M\prec\M$ is stationary.

We give our positive results first, which generalize Theorem 1.1 and Theorem 1.2 of \cite{Pillay-Yao-mini-flow}, in which $T$ is an $o$-minimal expansion of $\RCF$ and $H$ has dimension $\leq 1$.

\begin{theorem}\label{main-i}
Let $M$ be either $\Q$ or an $o$-minimal expansion of a real closed field, $G$ a group defined over $M$ admitting a ``$\dfg$-$\fsg$'' short exact sequence as in (\ref{equ-fsg-fsg-1}). Then
\begin{itemize}
    \item $H$ is stationary iff $G$ is stationary.
    \item (Global case) Suppose that either $G=H$ or $H$ is stationary, then every global  weakly generic type  is almost periodic.
    \item (Local case) Suppose that either $G=H$ or $H$ is stationary. Let $M_0=M^\ext$, then every weakly generic type  in $S_G(M_0)$ is almost periodic.
\end{itemize}
\end{theorem}

Let $\Gm$ be the multiplicative group of $\mathbb{M}\models \RCF$, then $\Gm$ and $\Gm^2$ are stationary $\dfg$ groups \cite{Petri}. If $\Gm$ is the multiplicative group of $\mathbb{M}\models pCF$, then each $\Gm^n$ is a stationary $\dfg$ group  for each $n\in \N^+$ \cite{Yao-f-generic-presb}.
Let $\Ga$ be the additive group of either $\mathbb R$ or $\Q$, then $\Ga$ is also a stationary  $\dfg$ group. By comparison, both $\Ga^2$ and  the borel subgroup $\Gm\rtimes \Ga$ of $\SL_2(\mathbb{M})$ are examples of non-stationary $\dfg$ groups. So it is reasonable to consider $\Ga^2$ and $\Gm\rtimes \Ga$ as ``minimal-non-stationary''  $\dfg$ groups.

Note that if $M$ is either
$(\Q,+,\times,0,1)$ or an $o$-minimal expansion of $(\mathbb R,+,\times,0,1)$,
then $M^\ext=M$ by \cite{Ded} and \cite{Delon}. In contrast to the above positive results, the next theorem gives  ``minimal counterexamples'', which generalizes  Theorem 1.3 of \cite{Pillay-Yao-mini-flow}.

\begin{theorem}
Let $M$ be $\mathbb R$ or $\Q$ in the language of rings. Let $H$ and $C$ be definable over $M$, where $H$ has $\dfg$ and $C$ has $\fsg$.
\begin{enumerate}
    \item If $H$ is bad (see Definition \ref{bad}) and $G=H\times C$, then working either in $S_G(M)$ or $S_G(\M)$ for $\M$ a monster model, the set of weakly generic types properly contains the set of almost periodic types.
    \item If $H$ is either $\Ga^2$ or a borel subgroup of $\SL_2(\mathbb{M})$, then $H$ is bad.
\end{enumerate}
\end{theorem}

The next  conjecture is based on the above result:

\begin{Conj}
Let $G$ be a definable, definably amenable group, defined  over either an $o$-minimal structure or a $p$-adically closed field. If $G$ is not $\dfg$, then the set of weakly generics coincides with the set of almost periodics iff $G$ is stationary.
\end{Conj}

\subsection{Notation and conventions}

$L$ will denote a language, $T$  a complete  theory, $M,N...$ models of $T$. It will be convenient to assume that $T$ has models which are $\kappa$-saturated and of cardinality $\kappa$ for arbitrarily large $\kappa$. Such a model, $\M$ say, will have homogeneity properties in addition to saturation properties. We call $\M$  a monster model. Let us fix a monster model $\M$ in this paper. A subset $A$ of $\M$ is called small if $|A|<|\M|$. We sometimes pass to an $|\M|^+$-saturated extension $\bar \M$ where types over $\M$ can be realized.  We usually write tuples as $a, b, x, y...$ rather than $\bar a, \bar b, \bar x, \bar y...$. A ``type'' is a complete type, and a ``partial type'' is a partial type. By a ``global type'' we mean a complete type over $\M$ (or $\bar \M$).  If $X$ is a definable set, defined over $M$ then we write $S_X(M)$ for the space of complete types concentrating on $X$. Let $\phi(x)$ be any $L_\M$-formula with $x=(x_1,..,x_n)$, and $A\sq \M$, then $\phi(A)$ is defined to be the set $\{a\in A^n|\ \M\models \phi(a)\}$. We sometimes use $X(x)$ to denote the formula which defines $X$, and identify $X$ with points in $\M$, namely $X=X(\M)$. If $A\sq B$ and $p\in S(B)$, then by $p|A$ we mean the  restriction of $p$ to $A$, namely, the collection of all formulas in $p$ with parameters from $A$. Suppose that $G$ is a group and $g,h\in G$, then by $g^h$ we mean the conjugate $h^{-1}gh$.

Our notation for model theory is standard, and we will assume familiarity with basic notions such as type spaces, heirs, coheirs, definable types etc. References are \cite{Pozit-Book} as well as \cite{Sim-Book}.

The paper is organized as follows:
For the rest of the introduction we give precise definition and preliminaries relevant to our results.

In Section 2,  we prove part (1) and (2) of Theorem \ref{main-i}, the positive results for global case.

In Section 3, we prove part (3) of Theorem \ref{main-i},  the positive result for Local case.

In Section 4, we show that when a $\dfg$ group $H$ is bad (see Definition \ref{bad}), then, for both global and local cases, we can find a weakly generic type on $H\times C$ which is not almost periodic whenever $C$ is an infinite $\fsg$ group. We also show that $\mathbb{G}_a^2$ and any borel subgroup of $\SL_2(\mathbb{M})$ are examples of bad $\dfg$ groups.

\subsection{Definable Topological dynamics}
Our reference for (abstract) topological dynamics are \cite{book-mini-flows-I,book-mini-flows-II}. Let $G$ be a topological (often discrete) group, by a \emph{$G$-flow} we mean  an  action $G\times X\rightarrow X$ of $G$ on  a  compact  Hausdorff  topological space $X$ by homeomorphisms, and denote it by $(G,X)$. Often it is assumed that there is a dense orbit.

By a \emph{subflow} of $X$ we mean a closed $G$-invariant subset of $X$. Minimal subflows of $X$ always exist.
A point $x\in X$ is \emph{almost periodic} if the closure $\cl(G\cdot x)$ of its $G$-orbit is  a  minimal subflow of $X$. Equivalently, $x\in X$ is almost periodic if $x$ is in some minimal subflow of $X$.

Given a $G$-flow $(G,X)$, its \emph{enveloping semigroup} $E(X)$ is the closure in the space $X^X$ (with the product topology) of the set of maps $\pi_g: X\rightarrow X$, where $\pi_g(x)=g\cdot x$, equipped with the composition (which is continuous on the left). So any $e\in E(X)$  is
 a map from $X$ to $X$.

\begin{fact}\label{Ellis-env}
Let $X$ be a $G$-flow. Then
\begin{itemize}
    \item $E(X)$ is also a $G$-flow and $E(E(X))\cong E(X)$ as $G$-flows.
    \item For any $x\in X$, the closure of its $G$-orbit is exactly $E(X)(x)$. Particularly, for any $f\in E(X)$, $E(X)\circ f$ is the closure of $G\cdot f$.
\end{itemize}
\end{fact}

In the model theoretic context,  we consider a complete theory $T$, model $M$ of $T$, group $G$ definable over $M$ and the action of $G(M)$ on the type-space $S_G(M)$ as $gp=\tp(ga/M)$ where $a$ realizes $p$. It is easy to see that $S_G(M)$ is a $G(M)$-flow with a dense orbit $\{\tp(g/M)|\ g\in G(M)\}$.

Take a monster model $\M$ and identify $G$ with $G(\M)$. We call a formula $\varphi(x)$,  with parameters in $\M$, a \emph{$G$-formula} if $\varphi(\M)$ is a definable subset of $G$. Suppose that $\varphi(x)$ is a $G$-formula   and $g\in G$, then  the left translate $g\varphi(X)$ is defined to be $\varphi(g^{-1}x)$.  It is easy to check that $(g\varphi)(\M)=gX$ if  $X=\varphi(\M)$. For $p\in S_G(M)$, we have $gp=\{g\varphi(x)|\ \varphi\in p\}$.

We recall some notions from \cite{Newelski-I}.
\begin{definition}
\begin{enumerate}
    \item A definable subset $X\sq G$ is generic if finitely many left translates of $X$ cover $G$. Namely, there are $g_1,...,g_n\in G$ such that $\bigcup_{i=1}^n g_iX=G$.
    \item A definable subset $X\sq G$
    is weakly  generic if there is an non-generic definable subset $Y$ such that $X\cup Y$ is generic.
    \item A formula $\varphi(x)$ is generic if the definable set $\varphi(\M)$ is generic (similarly for weakly generic formulas).
    \item A complete type $p\in S_G(A)$ is generic if every formula in $ p$ is generic.
    \item Likewise $p\in S_G(A)$ is weakly generic if every  formula in $ p$ is weakly generic.
\end{enumerate}
\end{definition}

\begin{fact}\cite{Newelski-I}\label{fact-New-min-flow}
\begin{itemize}
\item Let $\AP(S_G(M))\sq S_G(M)$ be the space of almost periodic types, and $\WG(S_G(M))\sq S_G(M)$ the space of weakly generic types. Then $\WG(S_G(M))=\cl(AP(S_G(M)))$.

\item If there is a generic type in $S_G(M)$, then there is a unique minimal subﬂow of $S_G(M)$ which moreover coincides with the set of generic types. So also generic types, almost periodic types, and weakly generic types coincide.
\end{itemize}
\end{fact}

Let $M^\ext$ be a  Shelah  expansion of $M$ in the language
$L^\ext=\{R_\phi(x)|\ \phi(x)\in L_\M\}$ with $R_\phi(M)=\phi(M)=\{\phi(a)|\ a\in M\}$, and $T_M^{\ext}$ the complete theory of $M^\ext$. We denote the collection of quantifier-free types over $M^\ext$ which concentrate on $G$ by $S_{G,\ext}(M)$. The space $S_{G,\ext}(M)$ is naturally homeomorphic to the space $S_{G,M}(\M)$ of global complete types  concentrating on $G$ which are finitely satisfiable in $M$, via the map
\[
p\in S_{G,M}(\M)\mapsto \{\psi(M)|\ \psi\in p \}\in S_{G,\ext}(M).
\]

\begin{fact}\label{Semigp-struc}\cite{Newelski-I}
The enveloping semigroup $E(S_{G,\ext}(M))$  of  $S_{G,\ext}(M)$ is isomorphic to $(S_{G,M}(\M),*)$ where $*$ is defined as following: for any $p,q\in S_{G,M}(\M)$, $p *q=tp(b\cdot c/\M)$ with $b$ realizes $p$ and $c$ realizes $q$, and $tp(b/\M,c)$ is finitely satisfiable in $M$.
\end{fact}
\begin{rmk}
It is easy to see from Fact \ref{Ellis-env} and Fact \ref{Semigp-struc} that for any $p\in S_{G,M}(\M)$,
\[
\cl(G(M)\cdot p)= S_{G,M}(\M) *p.
\]
\end{rmk}

\subsection{$\NIP$, Definable amenability, and connected component}

Recall that a theory $T$ is said to be \emph{$\NIP$} if for any indiscernible sequence $(b_i:\ i<\omega)$, formula $\phi(x,y)$ and $a\in \M$, there is an eventual truth-value of $\phi(a,b_i)$ as $i\rightarrow \infty$. Now we assume that $T$ has NIP throughout this paper.

Let $\phi(x,y)$ be a formula. Recall that a formula $\phi(x,b)$ \emph{divides} over a set $A$ if there is an infinite $A$-indiscernible sequence $(b=b_0,b_1,b_2,...)$ such that $\{\phi(x,b_i)|\ i<\omega\}$ is inconsistent. A type $p\in S(B)$ divides over $A\sq B$ if there is a formula $\phi\in p$ divides over $A$. By \cite{HP-$NIP$-inv-measure},
a global type $p\in S(\M)$ does not divide over a small model $M$ if and only if $p$ is $\Aut(\M/M)$-invariant.

\begin{fact}\label{fact-SG(M^ext)}\cite{HPP-$NIP$-JAMS}
 Let $M$ be a  model of $T$ and $G$  a definable group defined over $M$.
Then   $T^\ext_M$ has quantifier elimination and $\NIP$. So $S_{G,ext}(M)$ coincide with the space of $S_G(M^\ext)$ of complete types over $M^\ext$ concentrating on $G$.
\end{fact}

\begin{rmk}
We see from Fact \ref{Semigp-struc} and Fact \ref{fact-SG(M^ext)} that the semigroup operation ``$*$'' on $S_G(M^\ext)$ is  defined as following: for any $p,q\in S_G(M^\ext)$, $p *q=tp(b\cdot c/M^\ext)$ with $b$ realizes $p$ and $c$ realizes the unique heir of  $q$ over $M^\ext,b$.
\end{rmk}

Let $G=G(\M)$ be a definable group. Recall that a type-definable over $A$ subgroup $H$ is  a  type-definable over $A$ subset of $G$, which is also a subgroup of $G$. We say that $H$ has bounded index if $|G/H|<2^{|T|+|A|}$. For groups definable in  $\NIP$ structures, the smallest type-definable subgroup of bounded index exists \cite{HPP-$NIP$-JAMS}, which is the intersection of all type-definable subgroups of bounded index, we write it as $G^{00}$, and call it the \emph{type-definable connected component}. Another model theoretic invariant is $G^{0}$, called the \emph{definable-connected component} of $G$, which is the intersection of all definable subgroups of $G$ of   finite index.  Clearly, $G^{00}\leq G^0$.

\begin{fact}\cite{CPS-Ext} If  $G$ is definable over $M$, then $G^{00}$ is the same whether computed in $T$ or $T^{\ext}_M$.
\end{fact}

Recall also that a  \emph{Keisler measure} over $M$ on $X$, with $X$ a definable subset of $M^n$, is a finitely additive measure on the Boolean algebra of $M$-definable subsets of $X$. When we take the monster model, i.e. , $M=\M$,  we call it a global Keisler measure. A definable group $G$ is said to be \emph{definably amenable} if it admits a global (left) $G$-invariant probability  Keisler measure. By \cite{HPP-$NIP$-JAMS} this is equivalent to the existence of a $G(M)$-invariant probability  Keisler measure over $M$ on $G$, whenever $M$ is a model over which $G$ is defined.

A nice result of \cite{CS-amenable-$NIP$-group} shows that:
\begin{fact} Let $G$ be definable over $M$. Then
 $G$ is definably amenable iff  there exists $p\in S_G(\M)$ such that for every $g\in G=G(\M)$, $gp$ does not divide over $M$.
\end{fact}

Following the notation of \cite{CS-amenable-$NIP$-group}  we call a type $p$ as in the right hand side a (global) \emph{strongly $f$-generic} on $G$ over $M$. A global type  is strongly $f$-generic  if it is strongly  $f$-generic over some small model.

Given a definable subset $X$ of $G$, we say that $X$ is \emph{$f$-generic} if for some/any model $M$ over which $X$ is  defined and any $g\in G$, $gX$ does not divide over $M$. Call a complete type  $p$ (over some set of parameters) $f$-generic iff for every
formula $\psi(x)$ in $p$, $\psi(\M)$ is $f$-generic. In \cite{CS-amenable-$NIP$-group}, the authors showed that in $\NIP$ theories:
\begin{fact}
Let $G$ be a definable group. Then the following are equivalent:
\begin{enumerate}

    \item $G$ is definably amenable.
    \item $G$ admits a global type $p\in S_G(\M)$ with a bounded $G$-orbit.
    \item $G$ admits a global strongly $f$-generic type.

\end{enumerate}
\end{fact}

Moreover,
\begin{fact}\label{fact-st-f-generic} For a definably amenable group $G$, we have that
\begin{enumerate}
\item  Weakly generic definable subsets, formulas, and types coincide with $f$-generic definable subsets, formulas, and types, respectively.
\item $p\in S_G(\M)$ is  $f$-generic iff it has a bounded $G$-orbit.
\item $p\in S_G(\M)$ is $f$-generic if and only if it is $G^{00}$-invariant.
 A type-definable subgroup $H$ of the stabilizer of a global $f$-generic type is exactly $G^{00}$
 \item A global type is strongly $f$-generic over $M$ iff it is weakly generic and $M$-invariant, or  equivalently, does not divide over $M$.
\end{enumerate}
\end{fact}

Among the strongly $f$-generic types $p \in S_G(\M)$, there are two extreme cases:
\begin{enumerate}
\item There is a small submodel $M$ such that every left $G$-translate of $p$ is  finitely satisfiable in $M$, and we call such types the \emph{$\fsg$} (finitely satisfiable generic) types on $G$ over $M$;
\item There is a small submodel $M$ such that every left $G$-translate of $p$  is definable over $M$, and we  call such types the \emph{$\dfg$} (definable $f$-generic) types on $G$ over $M$.
\end{enumerate}

A definable group $G$ is called $\fsg$ or $\dfg$ if is has a $\fsg$ or $\dfg$ type, respectively.
Clearly, both $\fsg$ and $\dfg$ groups are definably amenable.
We now discuss these two cases. Let $\Stab_l(p)$ denote the stabilizer of $p\in S_G(M)$ with respect to the left group action, and $\Stab_r(p)$ the stabilizer of $p$ with respect
to the right group action. By \cite{HPP-$NIP$-JAMS} we have:

\begin{fact}\label{fact-fsg}
Let $G$ be an $\emptyset$-definable $\fsg$ group witnessed by a global type $p\in S_G(\M)$ and a small model $M$. Then
\begin{itemize}
    \item $p$ is both left and right generic.
    \item Any left (right) translate of $p$ is finitely satisfiable in any small model.
    \item $G^{00}=\Stab_l(p)=\Stab_r(p)$.
    \item Let $\Gen(G(\M))$ be the space global generic types in $S_G(\M)$, then $\Gen(G(\M))$ is the unique minimal subflow of $S_G(\M)$.
    \item Every global left/right generic type is $\fsg$.
    \item For any $N\prec\M$, every generic type $q\in S_G(N)$ has a unique global generic extension.
    \end{itemize}
    \end{fact}

\begin{rmk}\label{rmk-G00-inv-is-gen}
Let $G$ be an $\emptyset$-definable $\fsg$ group and $p\in S_G(\M)$ is $G^{00}$-invariant, then  by Fact \ref{fact-st-f-generic}  $p$ is weakly generic, thus is $p$ generic by Fact \ref{fact-New-min-flow}.
\end{rmk}
\begin{fact}\cite{CPS-Ext}\label{fact-fsg-ext}
Assume that $T$ is $\NIP$. Let  $G$ be an $\emptyset$-definable $\fsg$ group and $M\prec\M$.
    \begin{itemize}
    \item $G$ also has $\fsg$ when we compute it in $T^\ext_M$.

    \item  $q\mapsto \{\psi(M)|\ \psi\in q\}$ is a bijection between $\Gen(G(\M))$ and $\Gen(G(M^\ext))$.

    \item $\Gen(G(M^\ext))$ is a bi-ideal of (the semigroup) $S_G(M^\ext)$.
\end{itemize}
\end{fact}

We now discuss the $\dfg$ groups.

\begin{fact}\label{dfg-fact-1}\cite{Pillay-Yao-mini-flow}
Assume that $T$ is $\NIP$. Let   $G$ be an $\emptyset$-definable group and $p\in S_G(\M)$ an $f$-generic type. If $p$ is definable over  $M$, then
\begin{itemize}
    \item Every left  translate of $p$ is definable over $M$;
    \item $G^{00}=G^{0}=\Stab_l(p)$.
    \item $G\cdot p$ is closed, and hence a minimal subflow of $S_G(\M)$.
\end{itemize}
\end{fact}

\begin{fact}\label{dfg-fact-2}\cite{CPS-Ext}
Assume that $T$ is $\NIP$. Let $G$ be a $\dfg$ group definable over $M$. Let $N^*$ be an $|M|^+$-saturated extension of $M^\ext$. Then $p\in S_G(M^\ext)$ is almost periodic iff its unique  heir $\bar p$ over $N^*$ is weakly generic. Moreover,   any $G(N^*)$-translate of $\bar p$ is an heir of some $q\in S_G(M^\ext)*p$ over $N^*$. So $\bar p$ is a $\dfg$ type and thus $G$ is $\dfg$ with respect to $T^\ext_M$.
\end{fact}

\begin{fact}\label{min-flow-restriction}\cite{Sim-VC}
Assume that $T$ is $\NIP$, and $G$ is a definably amenable group definable over $M$. Let $M\prec N$, $\pi: S_G(N)\rightarrow S_G(M)$ the canonical restriction map, and $\cal M$  a minimal $G(N)$-subflow of $S_G(N)$. Then $\pi({\cal M})$ is a  minimal $G(M)$-subflow of $S_G(M)$.
\end{fact}

We conclude directly from Fact \ref{dfg-fact-2} and Fact \ref{min-flow-restriction}  that
\begin{coro}
Assume that $T$ is $\NIP$. Let $G$ be a $\dfg$ group definable over $M$, $N^*$   an $|M|^+$-saturated extension of $M^\ext$, and $\cal J$ a minimal $G(N^*)$-subflow of $S_G(N^*)$. Then $\pi({\cal J})$ is a minimal $G(M)$-subflow of $S_G(M^\ext)$ and $\pi$ is a bijection from $\cal J$ to $\pi({\cal J})$.
\end{coro}

The research is supported by The National Social Science Fund of China(Grant No. 20CZX050). The authors would like to thank the referee for carefully reading our paper and offering detailed comments, which have been very helpful for us as we revise the paper.

\section{Stationary $\dfg$ groups and the global case}

Here we prove  Theorem \ref{main-i}. We assume that $M$ is $\Q$ or an $o$-minimal expansion of a real closed field and $T=\Th(M)$. Assume that $\M$ be a monster model of $T$ and $G$ is an $M$-definable group defined in $\M$. We also assume that $G$ admits a $M$-definable short exact sequence
\begin{equation}\label{equ-sfg-fsg-2}
    1\rightarrow H\rightarrow G\rightarrow_\pi C\rightarrow 1,
\end{equation}
where $H$ has $\dfg$ and $C$ has $\fsg$.
As $T$ has definable Skolem functions (see \cite{definable Skolem functions} and Proposition 1.2 of Chapter 6 in \cite{L-van-den-Dries-Book}), let $f:C\rightarrow G$ be a definable section of $\pi$. We  aim to show that if $H$ is stationary, then every weakly generic type on $G$ is almost periodic.

\begin{rmk}\label{DA-by-DA}
  Note that in an arbitrary $\NIP$ theory, if $Y$ is a definable group and $X$ a definable normal subgroup of $Y$ such that both $X$ and $Y/X$ are definably amenable, then so is $Y$ (see  Exercise 8.23 of \cite{Sim-Book}). So any group $G$ admitting a decomposition as in (\ref{equ-sfg-fsg-2}) is definably amenable.
\end{rmk}

The following Lemma is formally analogous to Lemma 2.2 in \cite{Pillay-Yao-mini-flow}, replacing ``$H$-invariant'' with ``$f$-generic''.

\begin{lemma}\label{dfg-fsg-decomposition}
 Let $p=\tp(a/\M)$ be a weakly generic type on $G$ where $a\in G(\bar \M)$. Write $a$ uniquely as $h\cdot f(\pi(a))$ with $h\in H(\bar \M)$. Then $\tp(\pi(a)/\M)$ is a generic type on $C$ and $\tp(h/\M)$ is a weakly generic type on $H$.
\end{lemma}
\begin{proof}
Since $Stab(p)=G^{00}$, we see that $Stab(\pi(p))=\pi(G^{00})=C^{00}$ and so $\pi(p)$ is generic by Remark \ref{rmk-G00-inv-is-gen}.

Let $\eta: G(\bar \M)\rightarrow H(\bar \M)$ be the function given by $\eta(x)=x\cdot f(\pi(x))^{-1}$, then $h=\eta(a)$. For each $h_0\in H$, we have
\[
h_0\cdot h=h_0\cdot a\cdot f(\pi(a))^{-1}=(h_0\cdot a)\cdot f(\pi(h_0\cdot a))^{-1}=\eta(h_0\cdot a)\]

Since $\tp(a/\M)$ is weakly generic, we see that
\[
H\cdot \tp(a/\M)=\{\tp(h_0\cdot a/\M)|\ h_0\in H\}
\]
is bounded. So  the $H$-orbit
\[
\{\tp(h_0\cdot h/\M)|\ h_0\in H\}=\{\tp (\eta(h_0\cdot a)/\M)|\ h_0\in H\}
\]
of $\tp(h/\M)$ is bounded. So $\tp(h/\M)$ is a weakly generic as required.

\end{proof}

\begin{lemma}\label{every-heir-is-f-generic}
Let  $p\in S_H(M)$. If $p$ is a weakly generic type, then every global heir of $p$ is weakly generic.
\end{lemma}
\begin{proof}
If $H$ is defined in a model $M$ of an $o$-minimal expansion of $\RCF$, then $p$ is weakly generic iff it is has a global weakly generic extension $\bar p\in S_H(\M)$. By Fact \ref{fact-st-f-generic} and \ref{dfg-fact-1}, $\bar p$ is $H^{0}$-invariant, where $H^0$ is an $M$-definable group of $H$ by DCC (see \cite{Pillay-DCC}), so $p$ is $H^{0}(M)$-invariant, and thus every global heir of $p$ is also $H^0(\M)$-invariant, hence also  weakly generic. If $H$ is defined in an $p$-adically closed field, then by \cite{PY-dfg-groups} $H$ is eventually a trigonalizable algebraic groups. By Corollary 2.16 of \cite{Yao-tri}, we see that every  global heir of $p$ is weakly generic.
\end{proof}

\begin{fact}\label{fact-heir-AP}(\cite{Yao-tri},  Lemma 2.3)
Let $X\sq \M$ be a definably amenable group definable over $M\prec\mathbb{M}$. If every global heir of $p\in S_X(M)$
is  weakly generic, then $X(M)\cdot p$ is closed. In particular, $p$ is almost periodic.
\end{fact}

We conclude directly from Lemma \ref{every-heir-is-f-generic} and Fact \ref{fact-heir-AP} that

\begin{coro}
Let $p\in S_H(M)$. If $p$ is a weakly generic type, then   $p$ is almost periodic.
\end{coro}

\begin{lemma}\label{strongly-f-generic-extension}
Suppose that $p\in S_H(M)$ is weakly generic. Then $p$ has a global extension $  p^*\in S_H(\M)$ which is strongly $f$-generic over $M$.
\end{lemma}
\begin{proof}
Let $\bar p\in S_H(\M)$ be any weakly generic extension of $p$. Then $\bar p$ is almost periodic. By Proposition 3.31 of \cite{CS-amenable-$NIP$-group},
there is a global Keisler measure $\mu$ on $H$ such that $\mu(\phi(\M))>0$ for all $\phi\in \bar p$. By Lemma 5.8 of \cite{HP-$NIP$-inv-measure},
there is an $M$-invariant global Keisler measure $\mu^*$ on $H$ such that $\mu|M=\mu^*|M$. Let $p^*$ be an global extension of $p$ such that $\mu^*(\phi(\M))>0$ for each $\phi\in p^*$. So $\mu^*(h\phi(\M))>0$ for each $\phi\in p^*$ and $h\in H$. As $\mu^*$ is $M$-invariant, $\mu^*$ does not fork over $M$. We see that every left $H$-translate of $p^*$ does not fork over $M$. So $p^*$ is strongly $f$-generic over $M$ as required.
\end{proof}

Recall from the Section of Introduction that a definable group $X$ is stationary if every weakly generic type over any model has a unique global weakly generic extension. By Fact \ref{fact-fsg}, the $\fsg$ group $C$ is stationary.

In the rest of this section, we may assume that $G$ and its ``$\dfg$-$\fsg$'' short exact sequence are definable over $\emptyset$ by naming parameters.

\begin{lemma}\label{stationary}
The following are equivalent:
\begin{enumerate}
    \item $H$ is stationary.
    \item $H$ has boundedly many global weakly generic types.
    \item Every global weakly generic type is $\emptyset$-definable.
    \item $H$ has boundedly many strongly $f$-generic types.
\end{enumerate}
\end{lemma}
\begin{proof}
\begin{itemize}
    \item 1 $\Rightarrow$ 2: Let $N\prec\M$ be any small submodel. Assume that $H$ is stationary, then every global  weakly generic type $\bar p$ is the unique weakly generic extension of $\bar p|N$. As there are only boundedly many types over $N$, there are boundedly many global  weakly generic types.

    \item 2 $\Rightarrow$ 3: We first prove that every global weakly generic type is definable. Suppose that there are at most $\lambda<|\M|$ weakly generic types over $\M$. If  $p\in S_H(\M)$ is weakly generic but not definable, then $p$ has unboundedly many heirs  over an $|\M|^+$-saturated extension $\bar \M$ of $\M$ (see Proposition 1.19 in \cite{Pillay-Book-Stability}), and all of them are weakly generic by Lemma \ref{every-heir-is-f-generic}. Let $\{q_i|\ i<\lambda^+\}$ be a set of distinct weakly generic types over $\bar \M$. Let $\phi_{ij}(x,b_{ij})$ be the formula such that $\phi_{ij}(x,b_{ij})\in p_i$ and $\neg\phi_{ij}(x,b_{ij})\in p_j$. Let $N$ be a small submodel of cardinality $\lambda^{+}$ which contains the set  $\{b_{ij}| i\neq j\in\lambda^+\}$. We see that $p_i|N\neq p_j|N$ for all $i\neq j\in\lambda^+$. So there are at least $\lambda^+$ many weakly generic types over $N$. By the saturation of $\M$, take some $N'\prec \M$ such that $\tp(N)=\tp(N')$. So there are at least $\lambda^+$ many weakly generic types over $N'$. We conclude that there are   at least $\lambda^+$ many weakly generic types over $\M$, which is a contradiction.

    Now take any global weakly generic type $p\in S_H(\M)$ and any small submodel $N_0$.  If $p$ is not $\emptyset$-definable, then we see from the definability of $p$ that  $\{\sigma(p)|\ \sigma\in\Aut(\M)\}$ is unbounded. But each $\sigma(p)$ is weakly generic. So $p$ is $\emptyset$-definable.

     \item 3 $\Rightarrow$ 1: Let $N$ be a small submodel and $p\in S_H(M)$  weakly generic. Suppose that $q$ is a global weakly generic extension of $p$. Then $p$ and $q$ are definable over $N$. We see that $q$ is the unique heir of $p$. So the $q$ is the unique global weakly generic extension of $p$.

     \item It is easy to see that 2 $\Rightarrow$ 4. We now show that 4 $\Rightarrow$ 1.

     Suppose that $H$ is non-stationary. Let $N$ be any small submodel. Then there is a weakly generic type $p\in S_H(N)$ such that $p$ has unboundedly many global weakly generic extensions. For any cardinal $\lambda<|\M|$, take a sufficiently large  submodel $N'\succ N$ such that $|N'|=\lambda$ and $p$ has $\lambda$ many different weakly generic extensions $\{p_i|\ i<\lambda\}$ over $N'$. By Lemma \ref{strongly-f-generic-extension}, each $p_i$ has a strongly $f$-generic extension over $\M$. Hence there are unboundedly many strongly $f$-generic extensions of $p$.
\end{itemize}
\end{proof}

Simon \cite{Sim-distal} has isolated a notion, distality, meant to express the property that a $\NIP$ theory $T$ has “no stable part”, or is “purely unstable”. Examples include any o-minimal theory and $\pCF$ \cite{Sim-distal,Sim-fin}.

Let $p(x), q(y)$ be global types such that $p$ is definable over a small submodel and $q(y)$ is finitely satisfiable in a small submodel, then $p(x)\otimes q(y)= q(y)\otimes p(x)$ (see Lemma 2.23, \cite{Sim-Book}). By Lemma 2.16 in \cite{Sim-distal}, we have the following fact:
\begin{fact}\label{distal}
Assume $T$ is a distal $\NIP$ theory. Let $\M$ be  saturated, $N\prec \M$ a small submodel, $p(x)\in S(\M)$ definable over $N$, and $q(y)\in S(\M)$ finitely satisfiable in $N$. Then $p(x)$ and $q(y)$ are \emph{orthogonal}. Namely, $p(x)\cup q(y)$ implies a complete global type. In fact, if $a\models p$ and $b\models q$, then $\tp(a/\M,b)$ is the unique heir of $\tp(a/N)$ and $\tp(b/\M,a)$ is finitely satisfiable in $N$.
\end{fact}

We now consider the case where $H$ is stationary.
The following Lemma says that if $H$ is stationary, then we can exchange the positions of $h$ and $f(\pi(a))$ in the ``decomposition'' of $a=h\cdot f(\pi(a))$ as in Lemma \ref{dfg-fsg-decomposition}.

\begin{lemma}\label{fsg-dfg-decomposition}
 Let $p=\tp(a/\M)$ be a weakly generic type on $G$. Write $a$ uniquely as $f(\pi(a))\cdot h'$ with $h'\in H$. If $H$ is stationary, then $\tp(\pi(a)/\M)$ is a generic type on $C$ and $\tp(h'/\M)$ is a weakly generic type on $H$.
\end{lemma}
\begin{proof}
We see from   Lemma \ref{dfg-fsg-decomposition}  that   $\tp(a/\M)$  is of the form $\tp(h\cdot f(\pi(a))/\M)$ such that $\tp(h/\M)$ and $\tp(\pi(a)/\M)$ are global weakly generic types of $H$ and $C$ respectively. It suffices to show that $\tp(h^{f(\pi(a))}/\M)$ is a weakly generic type on $H$.

We now show that the $H$-orbit of $\tp(\tp(h^{f(\pi(a))}/\M))$ is bounded. By Lemma \ref{stationary} and Fact \ref{distal}, $\tp(h/\M)$ and $\tp(\pi(a)/\M)$  are orthogonal, we conclude that  $\tp(h/\M, \pi(a))$ is the unique heir of $\tp(h/\M)$.  So $\tp(h/\M, \pi(a))$ is a weakly generic type, and thus its $H(\dcl(\M,\pi(a)))$-orbit
\[
\{\tp(h_0h/\M, \pi(a))|\ h_0\in H(\dcl(\M,\pi(a)))\}
\]
is bounded. Since $x \mapsto  x^{f(\pi(a))}$ is a $\pi(a)$-definable automorphism of $H$, we see that $H(\dcl(\M,\pi(a)))$-orbit of $\tp({h}^{f(\pi(a))}/\M, \pi(a))$ is bounded.  So the $H$-orbit of $ \tp(h^{f(\pi(a))}/\M) $ is bounded as required.
\end{proof}

\begin{rmk}
Let $H$ be stationary. Then for any global weakly generic type $p_{_H}$ and $q_{_C}$ on $H$ and $C$ respectively, we can speak about the type $f(q_{_C})\cdot p_{_H}$ (resp. $p_{_H}\cdot f(q_{_C})$) which is defined to be $\tp(f(c^*)\cdot h^*/\M)$ (resp. $\tp(h^*\cdot f(c^*)/\M)$) where $c^*\models q_{_C}$ and $h^*$ realizes $p_{_H}$. By Lemma \ref{stationary} and Fact \ref{distal}, $p_{_H}$ and $q_{_C}$ are orthogonal, so $\tp(h^*/\M, c^*)$ is the unique heir of $p_{_H}$ over $\M,c^*$. Let $\Gen(C)$ be the space of generic types in $S_C(\M)$, and $Hp_{_H}$ the $H$-orbit of $p_{_H}$,  then $f(\Gen(C))\cdot (Hp_{_H})$ will denote the set
\[
\{f(q)\cdot p|\ q\in\Gen(C)\  \text{and}\  p\in Hp_{_H}\}.
\]
\end{rmk}

\begin{fact}\label{closure}\cite{Yao-tri}
Let $X$ be a definable group and  $p\in S_X(M)$. Then
\[
\cl(X\cdot p)=\{\tp(b\cdot c/M)|\ c\models p,\ b\in X,\   \text{and}\ \tp(b/M,c)\ \text{is finitely satisfiable in}\ M \}
\]
\end{fact}

\begin{lemma}\label{closure-of-the-orbit}
Let $r=f(q_{_C})\cdot p_{_H}$ be a weakly generic type on $G$, where $q_{_C}\in \Gen(C)$ and $p_{_H}\in S_H(\M)$ a definable $f$-generic type on $H$. Then
\[
\cl(G\cdot r)=f(\Gen(C))\cdot Hp_{_H}
\]
is a minimal $G$-flow, and in particular, $r$ is almost periodic.
\end{lemma}
\begin{proof}
It is easy to see that ``$\cl(G\cdot r)=f(\Gen(C))\cdot Hp_{_H}$'' implies $\cl(G\cdot r)$ is a minimal  subflow, since any $r'\in \cl(G\cdot r)$ is also of the form $r'=f(q')\cdot p'$ with $q'\in f(\Gen(C))$ and $p'\in Hp_{_H}=Hp'$ a definable $f$-generic type on $H$.

We now show that $\cl(G\cdot r)=f(\Gen(C))\cdot Hp_{_H}$.
Firstly, ``$\cl(G\cdot r)\sq \Gen(C)\cdot Hp_{_H}$'' is contained in the proof of Lemma 2.4 in \cite{Pillay-Yao-mini-flow}.

The new observation here is that $\Gen(C)\cdot Hp_{_H}\sq \cl(G\cdot r)$. Let $q'\in \Gen(C)$, $p'\in H\cdot p_{_H}$, and $r'=f(q')\cdot p'$.
Let $\bar \M$ be an $|\M|^+$-saturated extension of $\M$ and $a^*\in G(\bar \M)$ realize $r$.
By Fact \ref{closure}, it suffices to show that there is $g\in G(\bar \M)$ such that $r'=\tp(ga^*/\M)$ where $\tp(g/\M,a^*)$ is finitely satisfiable in $\M$.

Since $\Gen(C)$ is the unique minimal subflow of $S_C(\M)$, we see that
\[
\cl(C\cdot q')=\cl(C\cdot q_{_C})=\Gen(C).
\]
Take any $b\in G(\bar \M)$ such that $q'=\tp(\pi(b)\pi(a^*)/\M)$, where $\tp(b/\M,a^*)$ is finitely satisfiable in $\M$.

Let $b=f(\pi(b))h$ and $a^*=f(\pi(a^*))h^*$, then
\[
ba^*=f(\pi(b))hf(\pi(a^*))h^*=f(\pi(b))f(\pi(a^*))h^{f(\pi(a^*))}h^*.
\]
Let $h_0\in H(\bar \M)$ such that $f(\pi(b))f(\pi(a^*))=f(\pi(b)\pi(a^*))h_0$. Then $h_0\in\dcl(b,\pi(a^*))$.
Since both $\tp(b/\M, \pi(a^*),h^*)$  and $\tp(\pi(a^*)/\M,h^*)$ are  finitely satisfiable in $\M$, we see that $\tp(h^*/\M,\pi(a^*),b)$ is the unique heir of $p_{_H}$, which implies that
\[
\tp(h_0h^{f(\pi(a^*))}h^*/\M)\in \cl(H\cdot p_{_H})=H  p_{_H},
\]
By the orthogonality of $\tp(\pi(b)\pi(a^*)/\M)\in \Gen(C)$ and $\tp(h_0h^{f(\pi(a^*))}h^*/\M)\in H p_{_H}$, we have that
\[
\tp(ba^*/\M)=\tp(f(\pi(b)\pi(a^*))h_0h^{f(\pi(a^*))}h^*/\M) = f(q')\cdot p,
\]
with $p\in Hp_{_H}$.

Take some small submodel $N$ such that $H(N)$ meets all coset of $H^0$. Then $H(N)^g$ also meets all coset of $H^0$ for each $g\in G(\bar \M)$ as $H^0$ is normal in $G$,   Let $h_*\models p$, $c_*\models q'$, and $c_0\in C$ such that $\tp(c_*/N )=\tp(c_0/N)$.  By orthogonality of $\tp(h_*/\M)$ and $\tp(c_*/\M)$,  we have that $\tp(c_*,h_*/N )=\tp(c_0,h_*/NM )$. For each $h_1\in H(N)$, we see that $\tp({(h_1h_*) }^{f(c_0)}/N)$ is weakly generic, and both $\tp({(h_1h_*) }^{f(c_*)}/\M)$ and $\tp({(h_1h_*) }^{f(c_0)}/\M)$ are weakly generic extensions of $\tp({(h_1h_*) }^{f(c_0)}/N)$.
By the stationarity of $H$, we conclude that
\[
\tp({(h_1h_* )}^{f(c_*)}/\M)=\tp({(h_1h_*)}^{f(c_0)}/\M).
\]
Since both $H(N)$ and $H(N)^{f(c_0)}$ meet every coset of $H^0$, we have $$p'\in Hp_{_H}=H(N)p_{_H}=H(N)^{f(c_0)}p_{_H}.$$  So there exist $h_s,h_t\in H(N)$ such that
 \[
 \tp(h_s^{{f(c_0)}}h_*/\M)=\tp(h_th_*/\M)=p'.
 \]
 It follows that
\begin{equation*}
\begin{split}
h_s\tp(h_*^{f(c_*)^{-1}}/\M)&=h_s\tp(h_*^{f(c_0)^{-1}}/\M)=\tp(({h_s}^{f(c_0)}h_*)^{{f(c_0)}^{-1}}/\M)\\
&=\tp(({h_t}h_*)^{f(c_0)^{-1}}/\M)=\tp(({h_t}h_*)^{f(c_*)^{-1}}/\M).
\end{split}
\end{equation*}
Let $g=h_sb$, then $\tp(g/\M,a^*)$ is finitely satisfiable in $\M$. We have that
\begin{equation*}
\begin{split}
\tp(ga^*/\M)&=h_sf(q')\cdot p=\tp(h_sf(c_*)h_*/\M)\\
&=(h_s\tp(h_*^{f(c_*)^{-1}}/\M))\cdot \tp(f(c_*)/\M)\ \ \ (\text{by orthogonality})\\
&=\tp(({h_t}h_*)^{f(c_*)^{-1}}/\M)\cdot \tp(f(c_*)/\M)\ \ \ (\text{by orthogonality})\\
&=\tp(f(c_*)/\M)\cdot \tp({h_t}h_*/\M)\ \ \ (\text{by orthogonality})\\
&=f(q')\cdot p'=r'.
\end{split}
\end{equation*}
This completes the proof.
\end{proof}

If $H$ is stationary, then by Lemma \ref{fsg-dfg-decomposition} every weakly generic $r\in S_G(\M)$ is of the form $f(q_{_C})\cdot p_{_H}$ with $q_{_C}\in \Gen(C)$ and $p_{_H}\in S_H(\M)$ a definable $f$-generic type on $H$. We conclude directly from Lemma \ref{closure-of-the-orbit} that:
\begin{theorem}\label{AP=WG-global}
If $H$ is stationary, then every global weakly generic type on $G$ is almost periodic.
\end{theorem}

\begin{lemma}\label{construct-strongly-f-generic}
 Let  $\tp(c^*/\M )$ be a generic type on $C$ and $\tp(h^*/\bar \M )$ is a strongly $f$-generic type on $H$, over $N\prec\M$, where $\bar \M$ is $|\M|^+$-saturated. Then $\tp(f(c^*)h^*/\M)$ is a strongly $f$-generic type on $G$ over $N$.
\end{lemma}
\begin{proof}
 For any $c_0\in C$ and $h_0\in H$, we have that
\[
f(c_0)h_0\cdot \tp(f(c^*)h^*/\M)= \tp(f(c_0)f(c^*)h_0^{f(c^*)}h^*/\M).
\]
Let $h_1=f(c_0c^*)^{-1}f(c_0)f(c^*)$ and $h_2=h_1h_0^{f(c^*)}$, then $h_2\in \dcl(\M,c^*)$ and
\[
\tp(f(c_0)f(c^*)h_0^{f(c^*)}h^*/\M)=\tp(f(c_0c^*)h_2h^*/\M)
\]
Since $\tp(h^*/\bar \M)$  is strongly $f$-generic over $N$, we have that $\tp(h_2h^*/\M,c_0c^*)$ does not fork over $N$. We also have that $\tp(c_0c^*/\M)$ is a $\fsg$ type on $C$, and thus does not fork over $N$. We conclude that  $\tp(h_2h^*, c_0c^*/\M)\in S_{H\times C}(\M)$ does not fork over $N$, so $\tp(f(c_0c^*)\cdot h_2h^*/\M)\in S_G(\M)$ does not fork over $N$ as required.
\end{proof}

\begin{prop}\label{H-sta=G-sta}
$H$ is stationary iff $G$ is stationary.
\end{prop}
\begin{proof}
Suppose that $H$ is stationary.  Let $N$ be a small elementary submodel  of $\M$ and  $p(x)=\tp(f(c)h/N)\in S_G(N)$ be any weakly generic type and $\bar p=\tp(f(c^*)h^*/\M)$ a global weakly generic extension of $p$ where $c^*\in C(\bar \M)$ and $h^*\in H(\bar \M)$. Then by Lemma \ref{fsg-dfg-decomposition}, we have that $c^*$ realizes a global generic type on $C$ and $h^*$ realizes a global definable weakly generic type on $H$. Since both $C$ and $H$ are stationary, we see that $\tp(c^*/\M)$ is the unique global generic extension of $\tp(c/N)$ and $\tp(h^*/\M)$ is the unique global  weakly generic extension of $\tp(h/N)$. By the orthogonality of $\tp(c^*/\M)$ and $\tp(h^*/\M)$, $\tp(f(c^*)h^*/\M)$ is determined by $\tp(c^*/\M)\cup\tp(h^*/\M)$. We conclude that $\bar p=\tp(f(c^*)h^*/\M)$ is the unique global weakly generic extension of $p$. So $G$ is stationary as required.

Conversely. Suppose that $H$ is not stationary. It suffices to show that $G$ has unboundedly many global weakly generic types. By Lemma \ref{stationary}, we see that $H$ has unboundedly many strongly $f$-generic types. Let  $\tp(c^*/\M )$ be a generic type on $C$ with $c^*\in C(\bar\M)$. Let $\{\tp(h_i/ \M )|\ i<\lambda\}$ be a set of distinct strongly $f$-generic types on $H$, with each $\tp(h_i/ \M )$   strongly $f$-generic over $N_i\prec\M$.   By Lemma 3.11 in \cite{CS-amenable-$NIP$-group}, each $\tp(h_i/ \M )$ extends to a  type $\tp(h_i^*/\bar \M)$, which is also strongly $f$-generic over $N_i$.  Then each $\tp(f(c^*)h^*_i/\M)$ is a strongly $f$-generic type on $G$ by Lemma \ref{construct-strongly-f-generic}.

Let $\eta(x)=f(\pi(x))^{-1}x$. Then $\eta$ is an is a $\emptyset$-definable function. Since each $h^*_i$ is $\eta((f(c^*))h^*_i)$, we see that $\tp(f(c^*)h^*_i/\M)\neq \tp(f(c^*)h^*_j/\M)$ when $i\neq j$.
So there are at least $\lambda$ many strongly $f$-generic types of $G$ for each $\lambda<|\M|$. We conclude that there are unboundedly many global weakly generic types on $G$, and this completes the proof.
\end{proof}

\begin{theorem}\label{stationary-AP=WG}
Assume the assumptions of Theorem \ref{main-i} hold. If $G$ is stationary, then the space of global almost periodics coincides with the space of global weakly generics.
\end{theorem}

\begin{proof}
It is immediately from Theorem \ref{AP=WG-global} and Proposition \ref{H-sta=G-sta}.
\end{proof}

\section{Local case}

We assume here that $M$ is an arbitrary model of $T$, where $T$ is  $\pCF$ or an $o$-minimal expansion of $\RCF$, in the language $L$.

We denote $M^\ext$ by $M_0$, and $\Th(M^\ext)$ by $T^\ext_M$. Let $\M_0$ be a monster model of $T^\ext_M$, and  $\M$ the restriction of $\M_0$ to $L$. Note that $\M$ is also a monster of $T$ since $\M$ is also a saturated model of cardinality arbitrary large.
We assume again that $G$ is a definable group over $M$, and admits  a $M$-definable short exact sequence
\[
1\rightarrow H\rightarrow G\rightarrow_\pi C\rightarrow 1,
\]
with $C$ a $\fsg$ group, $H$ a $\dfg$ group, and $f:C\rightarrow G$ a definable section of $\pi$. So we can write any $g\in G$ uniquely as $f(\pi(g))h$ or $h'f(\pi(g))$. Note that by Fact \ref{fact-fsg-ext} and Fact \ref{dfg-fact-2}, $H$ and $C$  also have $\dfg$ and $\fsg$, respectively, when we compute them in $T_M^\ext$.

The following Facts appear in \cite{CPS-Ext}.
\begin{fact}\label{AP-over-M0}
$p\in S_H(M_0)$ is almost periodic iff its unique global heir is a weakly generic type.
\end{fact}

\begin{fact}\label{defble-type-ext}
Suppose $p(x)\in S_H(M)$ is definable. Then $p(x)$ implies a unique complete type $p^*(x)\in S_H(M_0)$. Moreover, if $\bar p$ is the unique heir of $p$ over $\M$, then $\bar p$ implies a unique complete type over $\M_0$, which is precisely the the unique heir of $p^*$.
\end{fact}

Given  $p\in S(\M_0)$ we define $p_L=\{\phi(x,b)\in p| \ \phi\in L, \ b\in \M\}$, which is the restriction of $p$ to the language $L$.
\begin{lemma}\label{ext-stationary}
Suppose that $H$ is stationary  with respect to $T$. Let $ \WG(S_H(\M_0))$ and $\WG(S_H(\M))$ be the space of weakly generic types of $S_H(\M_0)$ and $S_H(\M)$ respectively. Then $p\mapsto p_L$ is a bijection from $ \WG(S_H(\M_0))$ to $\WG(S_H(\M))$. Particularly,   $H$ is stationary  with respect to $T^\ext_M$.
\end{lemma}
\begin{proof}
If $p\in S_H(\M_0)$ is weakly generic, then $H(\M_0)p$ is bounded, so $H(\M)p_L$ is also bounded, we see that   $p_L\in S_H(\M)$ is also weakly generic type on $H$, hence is definable over $M$.  By Fact \ref{defble-type-ext}, $p|M_0$ is determined by $p_L|M$, where  $p|M_0$  and   $p_L|M$ are the restrictions of $p$ and  $p_L$ to $M_0$ and $M$ respectively. Using Fact \ref{defble-type-ext} again, we see that $p_L$ determines a complete type over $\M_0$, which is the unique heir of $p|M_0$. It follows that $p$ is the unique heir of $p|M_0$. As every   weakly generic type $p$ on $H$ over $\M_0$ is the unique heir of $p|M_0$, we conclude that $H$ is stationary with respect to $T^\ext_M$.
\end{proof}

\begin{rmk}
We see from the previous Lemma that $H$ is stationary  with respect to $T$ iff it is stationary  with respect to $T^\ext_M$.
\end{rmk}

\begin{lemma}\label{orth0-in-M0}
Suppose that $H$ is stationary (in the sense of $T$ or $T_M^\ext$), $p(x)\in S_H(\M_0)$ is a  weakly generic type, $q(y)\in S_C(\M_0)$ is a finitely satisfiable generic type.  Then  $p(x)$ and $q(y)$ are  orthogonal.
\end{lemma}
\begin{proof}
Let $p(x)=\tp(h^*/\M_0)$ and $q(x)=\tp(c^*/\M_0)$. It suffices to show that $\tp(h^*/\M_0,c^*)$ is the unique heir of $\tp(h^*/\M_0)$. By Fact \ref{distal}, we see that  $p_L$ and $q_L$ are orthogonal. So $\tp(h^*/\M,c^*)$ is the unique heir of $\tp(h^*/M)$. By Fact \ref{defble-type-ext}, we see that $\tp(h^*/\M_0,c^*)$ is the unique heir of $\tp(h^*/M_0)$ as required.
\end{proof}

\begin{coro}\label{ext-stationary-AP=WG}
Suppose that $H$ is stationary (in the sense of $T$ or $T_M^\ext$). Then $G$ is stationary  with respect to $T^\ext_M$ and every global weakly generic type in $S_G(\M_0)$ is almost periodic.
\end{coro}
\begin{proof}
By Fact \ref{fact-fsg}, Fact \ref{fact-fsg-ext}, Lemma \ref{ext-stationary}, and Lemma \ref{orth0-in-M0},  we see that every weakly generic type   $r\in S_G(\M_0)$ is of the form $f(q_{_C})\cdot p_{_H}$, where $q_{_C}\in S_C(\M_0)$ is generic, $p_{_H}\in S_H(\M_0)$ is definable $f$-generic over $M_0$, and $q_{_C},\ p_{_H}$ are orthogonal.
A similar argument as in the proof of Lemma \ref{closure-of-the-orbit}  shows  that
\[
\cl(G(\M_0)\cdot r)=G(\M_0)\cdot r=f(\Gen(C))\cdot (H(\M_0)\cdot p_{_H}).
\]
So $r$ is almost periodic as required.  Similarly, the proof of Proposition  \ref{H-sta=G-sta} shows that $G$ is stationary  with respect to $T^\ext_M$.
\end{proof}

\begin{theorem}
Suppose that $H$ is stationary (in the sense of $T$ or $T_M^\ext$). Then every weakly generic type in $S_G(M_0)$ is almost periodic (in the sense of $T^\ext_M$).
\end{theorem}
\begin{proof}
By Corollary 4.7 in \cite{Sim-VC},
the restriction of a global almost periodic type to any submodel is also almost periodic. We conclude that every type in $S_G(M_0)$ is almost periodic by Corollary \ref{ext-stationary-AP=WG}.
\end{proof}

Also, with Proposition \ref{H-sta=G-sta}, it proves the local case of Theorem \ref{main-i}.

\

For the rest of this section,  $H$ need not to be  stationary.
Let $\I$ be the space of generic types in $S_C(M_0)$, then $\I$ is the unique minimal subflow of $S_C(M_0)$, which is also a bi-ideal of the semigroup $(S_C(M_0),*)$. We use $f(\I)$ to denote the set $\{f(q)|\ q\in \I\}$. Let $\J$ be the union of all minimal subflow of $S_H(M_0)$. For any $p\in \J$, $\J(p)$ denotes the minimal subflow generated by $p$.
We are going to describe the space of
almost periodic types in $S_G(M_0)$ via $\I$ and $\J$.

\begin{lemma}\label{bi-ideal}
$\J$ is a bi-ideal of $S_H(M_0)$.
\end{lemma}
\begin{proof}
Clearly, $\J$ is a left ideal. We now show that $\J$ is also a right ideal. Let $p_0\in \J$ and $p\in S_H(M_0)$. Then it suffices to show that $p_0*p$ is almost periodic. It is easy to see that
\[
\cl(H(M_0)\cdot p_0*p)=S_H(M_0)*p_0*p=\J(p_0)*p.
\]
For any $p_1\in \J(p_0)$, we have that
\[
\cl(H(M_0)\cdot p_1*p)=S_H(M_0)*p_1*p=\J(p_1)*p=\J(p_0)*p.
\]
We conclude that $\cl(H(M_0)\cdot p_0*p)$ is a minimal $H(M_0)$-flow, and hence $p_0*p$ is almost periodic as required.
\end{proof}

\begin{lemma}\label{Min-flow=I-J}
 For any $p\in \J$, $f(\I)*\J(p)$ is a minimal subflow of $S_G(M_0)$.
\end{lemma}
\begin{proof}
It suffices  to show that $f(\I)*\J(p)\sq S_G(M_0)*f(r)*s$ for any $r\in \I$ and $s\in \J(p)$. Let $u\in \I$ and $v\in \J(p)$, then there are $c\in C$ and $h\in H$ such that $u=\tp(c/M_0)*r$ and $v=\tp(h/M_0)*s$.
We assume that $b$ realizes the heir of $r$ over $\dcl(M_0, c)$, $\tp(h/M_0,c,b)$ is the heir of $\tp(h/M_0)$, and $h^*$ realizes the global heir of $s$. Let $h'\in H$ such that $f(c)f(b)h'=f(cb)$.
Then
\begin{equation*}
\begin{split}
f(u)*v&= \tp(f(cb)hh^*/M_0)=\tp(f(c)f(b)h'hh^{*}/M_0)\\
&=\tp(f(c)f(b)/M_0)*\tp(h'hh^{*}/M_0) \ \ \ \ (\text{by Fact \ref{dfg-fact-2}})\\
&=\tp(f(c)/M_0)*\tp(f(b)/M_0)*\tp(h'hh^{*}/M_0)
\end{split}
\end{equation*}
Let  ${h_0}\in (hh')^{f(b)^{-1}}H^0$ such that $\tp(h_0/M_0,b,h^*)$ is finitely satisfiable in $M_0$. Let $\bar c$ realize a   coheir of $\tp(c/M_0)$ over $M_0,h_0,b,h^*$. Then
\begin{equation*}
\begin{split}
&\tp(f(\bar c)/M_0)*\tp(h_0/M_0)*\tp(f(b)/M_0)*\tp(h^{*}/M_0)\\
&=\tp(f(\bar c)h_0f(b) h^{*}/M_0)\\
&=\tp(f(\bar c)f(b) {h_0}^{f(b)} h^{*}/M_0)\\
&=\tp(f(\bar c)f(b)/M_0)*\tp( {h_0}^{f(b)} h^{*}/M_0)\\
&=\tp(f(\bar c)/M_0)*\tp(f(b)/M_0)*\tp( hh' h^{*}/M_0)\\
&=\tp(f(  c)/M_0)*\tp(f(b)/M_0)*\tp( hh' h^{*}/M_0)
\end{split}
\end{equation*}
So $f(u)*v=\tp(f(  c) h_0/M_0)*f(r)*s$ as required.
\end{proof}

\begin{lemma}\label{AP=q-p-r}
Let $r\in S_G(M_0)$. Then $r$ is almost periodic iff $r=f(q)*p*r$ for some generic type $q\in \I$ and almost periodic type $p\in \J$.
\end{lemma}
\begin{proof}
For any $r\in S_G(M_0)$ and $p\in\J$, the previous lemma shows that $\cl(G(M_0)\cdot f(\I)*p*r)=f(\I)*\J(p)*r$. It is also easy to see that $\cl(G(M_0)\cdot f(q')*p'*r)=f(\I)*\J(p')*r=f(\I)*\J(p)*r$ for any $q'\in \I$ and $p'\in\J(p)$. So we conclude that $f(\I)*\J(p)*r$ is a minimal subflow. As $\cl(G(M_0)\cdot r)\supset f(\I)*\J(p)*r$, it follows that $r$ is almost periodic iff $r\in f(\I)*\J(p)*r$, which completes the proof.
\end{proof}

We now consider the case where $G$ is a product of $C$ and $H$. We identify $C$ with $C\times \{1_H\}$ and $H$ with $\{1_C\}\times H$, subgroups of $G=C\times H$

\begin{lemma}\label{prod-r=qp}
Suppose that $G=C\times H$. Then $r\in S_G(M_0)$ is almost periodic iff $r\in \I*\J$.
\end{lemma}
\begin{proof}
Let $q\in S_C(M_0)$ and $p\in S_H(M_0)$. We see from Lemma \ref{Min-flow=I-J} that $q*p$ is almost periodic if $q\in \I$ and $p\in \J$.

Conversely, suppose that $r\in S_G(M_0)$ is almost periodic, then by Lemma \ref{AP=q-p-r} we have $r=q*p*r$ for some $q\in \I$ and $p\in \J$. Let $N\succ M_0$ be any $|M_0|^+$-saturated extension. Take $c\in C(N)$ and $h\in H(N)$ such that $ ch$ realizes $r$. Let $h^*\in H(N)$ realize the coheir of  $p$ over $M_0,c,h$ and $c^*\in C$ realize the coheir of $q$ over $N$. Then
\[
q*p*r=\tp(c^*h^*ch/M_0)=\tp(c^*ch^*h/M_0),
\]
which is precisely $\tp(c^*c/M_0)*\tp(h^*h/M_0)$ since $\tp(c^*c/N)$ is a generic type over $N$. It is easy to see from Lemma \ref{bi-ideal} that $\tp(h^*h/M_0)=\tp(h^*/M_0)*\tp(h/M_0)$ is  almost periodic in $H$. Hence $r\in \I*\J$ as required.
\end{proof}

\section{Bad $\dfg$ groups and examples}
In this section, we assume that $M$ is the field $\mathbb R$ of real numbers, or $\Q$ of $p$-adic numbers, and $\M\succ M$ a monster model. As every type over $M$ is definable, we have that $M^\ext=M$.

\begin{definition}\label{bad}
Let $H$ be a $\dfg$ group definable over $M$. We say that $H$ is bad if there are  $\tp(a/\M)$ strongly $f$-generic over $M$ and an $M$-definable function $\theta$ such that $\tp(\theta(a)/\M)$ is a non-realized type finitely satisfiable in $M$.
\end{definition}

\begin{rmk}
The badness of $H$ is witnessed by any monster model $\bar \M\succ \M$.  Let $\tp(a/\M)$ be as in Definition \ref{bad}, then $\tp(a/\M)$ is $M$-invariant by Fact \ref{fact-st-f-generic} (4). By Lemma 3.11 of \cite{CS-amenable-$NIP$-group}, $\tp(a/\M)$ has an extension $\tp(a^*/\bar \M)$  which is also
strongly $f$-generic over $M$, so is $M$-invariant too. We conclude that $\tp(\theta(a^*)/\bar \M)$ is an $M$-invariant extension of $\tp(\theta(a)/\M)$. By the saturation of $\M$,   $\tp(\theta(a^*)/\bar \M)$  is also finitely satisfiable in $M$.
\end{rmk}
 We aim to show that every bad $\dfg$ group $H$ yields a counterexample $G=H\times C$, where there is weakly generic type on $G$ which is not almost periodic whenever $C$ is an infinite $\fsg$ group. We will also give two ``minimal'' examples of non-stationary $\dfg$ groups, and show that they are bad.

We denote the additive group and the multiplicative group of $\M$ by $\Ga$ and $\Gm$, respectively.
We use $|a|$ to denote the norm (or absolute value) of $a\in \M$. Note that relation $|x|\leq |y|$ is definable in both $\Th(\mathbb R)$ and $\Th(\Q)$ (see \cite{Mac}) in the language of rings. For any $a\in \M$, we say that  $a$ is \emph{bounded} over $M$ if $|a|<|b|$ for some $b\in M$. If $a\in \M$ is bounded over $M$, then there is $\st(a)\in M$ which is infinitesimally close to $a$ over $M$, namely, $0\leq |a-\st(a)|<|b|$ for all $b\in M\backslash\{0\}$. We call $\st(a)$ the \emph{standard part} of $a$. Clearly,  $a\in \M$ is unbounded over $M$ iff $a^{-1}$ is bounded and $\st(a^{-1})=0$.

The following  is a folklore. Nevertheless we give a proof here for convenience.
\begin{fact}\label{heir-neq-coheir}
Let $e\in \M\backslash M$  and $p\in S_1(M)$ be a non-algebraic type. Suppose that $p_1$ is the unique heir  $p$ over $M,e$, then $p_1$ is not a finitely satisfiable in $M$.
\end{fact}
\begin{proof}
Suppose that $p(x)=\tp(a/M)$. Then $a$ is either unbounded over $M$ or infinitesimally close to $\st(a)$ over $M$. Suppose for example that $a$ is  unbounded over $M$.    Let   $a^*\models p_1(x)$, then $a^*$ is also unbounded over $\dcl(M,e)$ since $p_1$ is the heir of $p$.   Since $e\notin M$, we see that either $e$ or $d=(e-\st(e))^{-1}$ is unbounded over $M$. Now $|a^*|>\max\{|e|,|d|\}$ but the formula ``$|x|>\max\{|e|,|d|\}$'' is not satisfiable in $M$. So $\tp(a^*/M,e)$ is not finitely satisfiable in $M$.
\end{proof}

\begin{prop}
Let $H$ be a bad $\dfg$ group, $C$ a definably compact group over $M$ with $\dim(C)\geq 1$, and $G=C\times H$. Then $S_G(M)$ has a weakly generic type which is not almost periodic.
\end{prop}
\begin{proof}
Let $\bar \M\succ\M$ be  $|\M|^+$-saturated. Let $c^*\in  C(\bar \M)$ realize a generic type on $C$ over $\M$. Let $\tp(h^*/\bar \M)\in S_H(\bar \M)$ be a strongly $f$-generic type on $H$ over $M$ such that $\theta(h^*)$ is finitely satisfiable in $M$, where $\theta$ is an $M$-definable function. We see from Lemma \ref{construct-strongly-f-generic} that $\tp((c^*,h^*)/\M)$ is strongly $f$-generic over $M$. By Fact \ref{fact-st-f-generic}, $\tp((c^*,h^*)/\M)$ is weakly generic.

Suppose for a contradiction that $\tp((c^*,h^*)/M)$ is almost periodic. Then by Lemma \ref{prod-r=qp}, $\tp(c^*/M,h^*)$ is finitely satisfiable in $M$. We conclude that both $\tp(c^*/M,\theta(h^*))$ and $\tp(\theta(h^*)/M,c^*)$ are finitely satisfiable in $M$, which contradicts to Fact \ref{heir-neq-coheir}.
\end{proof}

We say that $\tp(a/M)\in S_1(M)$ is a \emph{type of infinite} (resp. \emph{type of infinitesimal}) if and  $|a|>|b|$  for all $b\in M$ (resp. $0<|a|<|b|$  for all nonzero $b\in M$).

We now consider a definably compact subgroup $D$ of $M$, where \[
D=\SO_2(M)=\{(x,y)|\ x^2+y^2=1\}
\]
if $M=\mathbb R$, and
\[
D=\Z_p=\{x|\ |x|\leq 1\}
\]
if $M=\Q$. We tend to identify a point of $\SO_2$  with its $x$-coordinate,  working above the $x$-axis.

The following Facts can be found in \cite{PPY-SL2Qp} and \cite{YL}:

\begin{fact}
\begin{itemize}
    \item $D$ has $\fsg$.
    \item Let $p\in S_1(M)$ be infinitesimally close to $0\in M$. Then its global coheir  is a generic type on $D$, and its global heir is a definable $f$-generic type on $\Gm$.
\end{itemize}

\end{fact}

\begin{fact}
If  $a\in\bar \M$ is unbounded over $\M$ ($|a|>|b|$ for all $b\in\M$). Then $\tp(a/\M)$ is a definable $f$-generic type over $\emptyset$, on both $\Gm$ and $\Ga$.
\end{fact}

\begin{fact}
$\Gm$ and $\Ga$ are stationary $\dfg$ groups.
\end{fact}

\begin{fact}
If $D=\SO_2$, then
\[
D^{00}=\{x\in \M| \ x\ \text{is infinitesimal close to}\ 1 \ \text{over}\ M\}.
\]
If $D=\Z_p$, then
\[
D^{00}=\{x\in \M| \ x\ \text{is infinitesimal close to}\ 0\ \text{over}\ M \}.
\]
\end{fact}

A result of \cite{JY-abelian-pcf-group} shows that
\begin{fact}\label{fact-dfg-extension}
(Assuming $\NIP$) Suppose that $Y$ is a definable group and $X$ is a definable normal subgroup of $Y$. If both $X$ and $Y/X$ have $\dfg$, then $Y$ has $\dfg$.
\end{fact}

\begin{lemma}\label{f-generic-Ga2}
${\Ga\times \Ga}$ is a bad $\dfg$ group.
\end{lemma}
\begin{proof}
Clearly, $\Ga\times\Ga$ has $\dfg$ by Fact \ref{fact-dfg-extension}.
Let $p\in S_1(M)$ be a type of infinite, $a^*$ realize the global coheir of $p$, $b^*$ realize the global heir  of  $p$, and $h^*=a^*b^*$. Then it sufficies to show that $\tp(b^*, h^*/\M)\in S_{\Ga\times \Ga}(\M)$ is a strongly $f$-generic type over $M$.

Since both $\tp(a^*/\M)$ and $\tp(b^*/\M)$ are $M$-invariant, so is $\tp(a^*,b^*/\M)$ by the orthogonality of $\tp(a^*/\M)$ and $\tp(b^*/\M)$. Now $(b^*,h^*)$ and $(a^*,b^*)$ are interdefinable, we conclude that $\tp(b^*, h^*/\M)$ is also $M$-invariant.

We now prove that $\tp(b^*, h^*/\M)$ is weakly generic. Note that ${\Ga^2}={(\Ga^2)}^{00}$, so it suffices to show that
\[
\tp(a^*,b^*/\M)=\tp(\frac{h+h^*}{b+b^*},b+b^*/\M)
\]
for any $b,h\in \M$. By the orthogonality, we only need to prove that $\tp(a^*/\M)=\tp(\frac{h+h^*}{b+b^*}/\M)$. Clearly,  both $\tp(b^*/\M)$ and $\tp(h^*/\M)$ are unbounded over $\M$, and hence  definable $f$-generic types  on $\Ga$.

It is easy to see that $\tp((b/h^*)/\M)$ is infinitesimally close to $0$ over $\M$, thus is  a definable $f$-generic type  on $\Gm$. We now consider $\tp({a^*}^{-1}/\M)$ as global generic type on $D$. The it is easy to see  from the orthogonality that $\tp({a^*}^{-1}/\M,(b/h^*))$ is a generic extension of $\tp({a^*}^{-1}/\M)$, thus is finitely satisfiable in $M$, and $\tp((b/h^*)/\M, a^*)$ is   a definable $f$-generic type on $\Gm$ (definable over $M$).

Keep in mind that $\tp({a^*}^{-1}/\M)$ is a global generic type on $D$.

\vspace{0.2cm}

\noindent $\bullet$ If $D$ is $\SO_2$, then ${a^*}^{-1}+(b/h^*)={a^*}^{-1}(1+(a^*b/h^*))$. It is easy to see from the orthogonality that $\tp((a^*b/h^*)/\M,a^*)$  is  a definable $f$-generic type  on $\Gm$  over $M$, which is also infinitesimally close to $0$ over $\dcl(\M,a^*)$. We conclude that $1+(a^*b/h^*)\in D^{00}$ and thus
\[
\tp({a^*}^{-1}/\M)=\tp({a^*}^{-1}(1+\frac{a^*b}{h^*})/\M)=\tp({a^*}^{-1}+\frac{b}{h^*}/\M)
\]
as $\tp({a^*}^{-1}/\M, (a^*b/h^*))$ is generic.

\vspace{0.2cm}

\noindent $\bullet$ If $D$ is $\Z_p$, then $\tp({a^*}^{-1}/\M)=\tp({a^*}^{-1}+(b/h^*)/\M)$ since $b/h^*\in D^{00}$ and $\tp({a^*}^{-1}/\M, (b/h^*))$ is   generic.

\vspace{0.2cm}
So we conclude that
\[
\tp({a^*}^{-1}/\M)=\tp({a^*}^{-1}+\frac{b}{h^*}/\M)=\tp(\frac{b+b^*}{h^*}/\M).
\]
On the other side,
\[
\frac{h+h^*}{b+b^*}=\frac{h^*}{b+b^*}+\frac{h}{b+b^*}=\frac{h^*}{b+b^*}(1+\frac{b+b^*}{h^*}\frac{h}{b+b^*}).
\]
Using orthogonality again, we see that
$\tp(1+\frac{b+b^*}{h^*}\frac{h}{b+b^*}/\M)$ is infinitesimally close to $1$ over $\M$. So  $\tp((1+\frac{b+b^*}{h^*}\frac{h}{b+b^*})^{-1}/\M)$ is also infinitesimally close to $1$ over $\M$. We conclude that
\[
\tp(\bigg(\frac{h^*}{b+b^*}(1+\frac{b+b^*}{h^*}\frac{h}{b+b^*})\bigg)^{-1}/\M)=\tp(\frac{b+b^*}{h^*}(1+\frac{b+b^*}{h^*}\frac{h}{b+b^*})^{-1}/\M)=\tp({a^*}^{-1}/\M),
\]
which implies that
\[
\tp(a^*/\M)=\tp(\frac{h^*}{b+b^*}(1+\frac{b+b^*}{h^*}\frac{h}{b+b^*})/\M)=\tp(\frac{h+h^*}{b+b^*}/\M).
\]
We have that $\tp(b^*, h^*/\M)$ is weakly generic and $M$-invariant, so it is strongly $f$-generic over $M$ as required.
\end{proof}

Let $H$ be the borel subgroup of $\SL_2(\mathbb{M})$  consisting of the matrices of the form
$\begin{pmatrix}
t & u\\ 0 & t^{-1}
\end{pmatrix}$. We identify
$\begin{pmatrix} t & u\\ 0 & t^{-1}\end{pmatrix}$ with a pair $(t,u)$. Then
the group operation is given by
\[
(t,u)\cdot (t',u')=(tt', tu'+{t'}^{-1}u).
\]

\begin{lemma}\label{f-generic-Bor}
Let $H$ be the borel subgroup of $\SL_2(\mathbb{M})$ consisting of the matrices of the form
$\begin{pmatrix}
t & u\\ 0 & t^{-1}
\end{pmatrix}$. Then $H$ is a bad $\dfg$ group.

\end{lemma}
\begin{proof}
Since $H$ is an extension of $\Ga$ by $\Gm$, it has $\dfg$ by Fact \ref{fact-dfg-extension}. Let $p\in S_1(M)$ be a type of infinite, $a^*$  realize the global coheir of $p$, and $t^*$ realize the global heir  of $p$. Let $u^*=a^*t^*$. Then it sufficies to show that $\tp(t^*, u^*/\M)\in S_{H}(\M)$ is a strongly $f$-generic type over $M$. As we have showed in Lemma \ref{f-generic-Ga2}, we only need to prove that
\[
\tp(a^*/\M)=\tp(\frac{tu^*+u({t^*}^{-1})}{tt^*}/\M)
\]
for any $ u\in \Ga$ and  $t\in \Gm^{0}$.

Since $\tp({t^*}^2/\M)$ is a global definable $f$-generic type on $\Gm$, we see that $\tp((t/u){t^*}^2/\M)$ is also a global definable $f$-generic type on $\Gm$. Now
\[
\frac{tu^*+u({t^*}^{-1})}{tt^*}=a^*+\frac{ u({t^*}^{-1})}{tt^*}=a^*(1+{a^{*}}^{-1}\frac{ u({t^*}^{-1})}{tt^*}).
\]
By orthogonality, letting $\epsilon=({a^{*}}^{-1}  u({t^*}^{-1}))/(tt^*)$, we have that $\tp(1+\epsilon/\M)$ is infinitesimally close to $1$ over $\M$, so is $\emptyset$-definable. Clearly,
\[
\tp(a^*/\M)=\tp(a^*(1+\epsilon)/\M)\iff\tp({a^*}^{-1}/\M)=\tp({a^*}^{-1}(1+\epsilon)^{-1}/\M).
\]
Since $\tp({a^*}^{-1}/\M)$ is a generic type on $D$ and $\tp((1+\epsilon)^{-1}/\M)\vdash D^{00}$, we see from the orthogonality that  $\tp({a^*}^{-1}/\M)=\tp({a^*}^{-1}(1+\epsilon)^{-1}/\M)$.
So
\[
\tp(a^*/\M)=\tp(a^*+\frac{ u({t^*}^{-1})}{tt^*}/\M)=\tp(\frac{tu^*+u({t^*}^{-1})}{tt^*}/\M),
\]
which completes the proof.
\end{proof}

 Finally, we conjecture that
\begin{Conj}
Assume that $H$ is a $\dfg$ group definable in an $o$-minimal structure or a $p$-adically closed field, then $H$ is bad iff $H$ is non-stationary.
\end{Conj}

\end{document}